\documentclass[a4paper,11pt,reqno]{amsart}

\usepackage{latexsym,dsfont,amssymb,amsmath,amsthm,amscd}

\usepackage[all]{xy}

\chardef\bslash=`\\ 

\numberwithin{equation}{section}

\newtheorem{theorem}{Theorem}[section]
\newtheorem{corollary}[theorem]{Corollary}
\newtheorem{lemma}[theorem]{Lemma}
\newtheorem{proposition}[theorem]{Proposition}
\theoremstyle{remark}

\theoremstyle{definition}

\newcommand\bp{\begin{proof}}
\newcommand\ep{\end{proof}}

\newcommand\Dhat{{\hat\Delta}}

\newcommand\CC{{\mathcal C}}
\newcommand\D{{\mathcal D}}
\newcommand\E{{\mathcal E}}
\newcommand\F{{\mathcal F}}

\newcommand\RR{{\mathcal R}}
\newcommand\U{{\mathcal U}}

\newcommand\g{{\mathfrak g}}
\newcommand\h{{\mathfrak h}}

\newcommand\Clg{\operatorname{Cl}({\mathfrak g})}

\newcommand\GL{\operatorname{GL}}
\newcommand\End{\operatorname{End}}
\newcommand\Hom{\operatorname{Hom}}

\newcommand\Mod{\operatorname{-Mod}_f}
\newcommand\Nat{\operatorname{Nat}}
\newcommand\SU{\operatorname{SU}}

\newcommand\tr{\operatorname{tr}}
\newcommand\Vect{\mathcal Vec}

\newcommand{\ad}{\operatorname{ad}}
\newcommand{\add}{{\widetilde\ad}}

\newcommand{\C}{{\mathbb C}}

\newcommand{\N}{{\mathbb N}}

\newcommand{\R}{{\mathbb R}}
\newcommand\Sp{{\mathbb S}}
\newcommand\T{{\mathbb T}}
\newcommand\Z{{\mathbb Z}}

\newcommand\eps{\varepsilon}

\newcommand\enu[1]{\smallskip\newline\makebox[5mm][l]{\rm(#1)}}

\begin{document}

\title[Symmetric invariant cocycles]
{Symmetric invariant cocycles on the duals of $q$-deformations}


\author[S. Neshveyev]{Sergey Neshveyev}
\address{Department of Mathematics, University of Oslo,
P.O. Box 1053 Blindern, NO-0316 Oslo, Norway and \\ Institut de
Math\'ematiques de Jussieu, Universit\'e Paris 7 Denis Dedirot
}

\email{sergeyn@math.uio.no}

\author[L. Tuset]{Lars Tuset}
\address{Faculty of Engineering, Oslo University College,
P.O. Box 4 St.~Olavs plass, NO-0130 Oslo, Norway}
\email{Lars.Tuset@iu.hio.no}

\thanks{Supported by the Research Council of Norway.}

\date{February 13, 2009; minor changes January 8, 2011}

\begin{abstract}
We prove that for $q\in\C^*$ not a nontrivial root of unity any
symmetric invariant $2$-cocycle for a completion of $U_q\g$ is the
coboundary of a central element. Equivalently, a Drinfeld twist
relating the coproducts on completions of $U_q\g$ and $U\g$ is unique
up to coboundary of a central element. As an application we show that
the spectral triple we defined in an earlier paper for the
$q$-deformation of a simply connected semisimple compact Lie group $G$ does
not depend on any choices up to unitary equivalence.
\end{abstract}

\maketitle

\bigskip

\section*{Introduction}

Let $G_q$, $q>0$, be the compact quantum group which is the
$q$-deformation in the sense of Drinfeld and Jimbo of a simply connected semisimple compact Lie group $G$. In \cite{NT2} we  constructed a quantum
Dirac operator $D_q$ on $G_q$ that defines a biequivariant spectral
triple, which is an isospectral deformation of that defined by the
Dirac operator $D$ on $G$. To do this we used an analytic version of a
result by Drinfeld, due to Kazhdan and Lusztig. This result, see
\cite{NT3} and references therein, produces what we call a Drinfeld
twist $\F$, which is an element in the group von Neumann algebra $W^*
(G\times G)$, and an isomorphism $\varphi\colon W^*(G_q )\to W^* (G)$
such that $(\varphi\otimes\varphi)\Dhat_q
=\F\Dhat\varphi(\cdot)\F^{-1}$ and
$(\varphi\otimes\varphi)(\RR)=\F_{21}q^t\F^{-1}$, where $\RR$ is the
universal $R$-matrix for $U_q\g$ and $t\in\g\otimes\g$ is defined by a
suitably normalized $\ad$-invariant symmetric form on $\g$. Most
importantly, the Drinfeld associator~$\Phi_{KZ}$, which is defined via
monodromy of the KZ-equations,  is the coboundary of $\F^{-1}$ with
respect to the coproduct on $U_q\g$.

From the outset $D_q$ and the associated spectral triple depend on the
choice of $(\varphi, \F)$. In this paper we show that a different
choice in fact produces the same spectral triple up to unitary
equivalence, see Theorem~\ref{coruniqueD}. So our construction is as
canonical as one could possibly hope for. Everything hinges on a
uniqueness result for the Drinfeld twist which we establish in
Theorem~\ref{corunique}. It states that for a fixed $\varphi$, any
Drinfeld twist has to be of the form $(c\otimes c)\F\Dhat(c)^{-1}$,
where $c$ is a unitary central element of $W^*(G)$. Combining this with
the contribution from choosing a different $\varphi$, we deduce that
any Drinfeld twist is of the form $(u\otimes u)\F\Dhat(u)^{-1}$ for a
unitary $u\in W^*(G)$.

An equivalent form, Theorem~\ref{thmMain}, of the uniqueness result for
the Drinfeld twist says that any unitary symmetric $G_q$-invariant
$2$-cocycle in $W^*(G_q\times G_q)$, see Section~\ref{seccohom} for
precise definitions, is the coboundary of a central element. The result
is also true for nonunitary cocycles, and in this form it makes sense
and is valid for all $q\in\C^*$ not a nontrivial root of unity.

\smallskip

The study of $2$-cocycles on duals of compact groups was initiated by
Landstad~\cite{La} and Wassermann~\cite{Wa1}. They showed that
cohomology classes of unitary $2$-cocycles are in a one-to-one correspondence
with full multiplicity ergodic actions on operator algebras. It is
expected~\cite{Wa2} that any unitary $2$-cocycle on~$\hat G$ sufficiently close
to the trivial one is defined by a $2$-cocycle on the dual of a maximal
torus, but this has been proved only for some low rank groups. In
particular, it has been shown that $H^2(\widehat{\SU(2)};\T)$ is trivial.
The theory of full multiplicity ergodic actions was extended to compact
quantum groups in~\cite{BDRV}, and the second cohomology was computed
for the duals of free orthogonal quantum groups, which implies that
$H^2(\widehat{\SU_q(2)};\T)$ is trivial. Therefore for $\SU_q(2)$ we already know that any unitary $2$-cocycle is a coboundary, but our result says that if the cocycle is in addition invariant and symmetric, then it is a coboundary of a central element.

\smallskip

The paper is organized as follows. After a brief introduction to
nonabelian cohomology for Hopf algebras, we state our main result,
Theorem~\ref{thmMain}. In Section~\ref{smain} we show that any
symmetric invariant $2$-cocycle is cohomologous to a cocycle $\E$
satisfying two additional properties; that it acts trivially on the
isotypic components of the tensor product of two modules corresponding
to the highest and next to highest weights. Our goal then is to show
that $\E=1$. This can easily be done for $SU_q(2)$ because the fusion
rules are sufficiently simple. However, for higher rank groups new
impetus is needed.

Our approach is motivated by the work of Kazhdan and
Lusztig~\cite{KL1}, who constructed a comonoid in a completion of the
Drinfeld category representing the forgetful functor. Such a comonoid,
in the equivalent category of $U_q\g$-modules, is also implicit in
Lusztig's book~\cite{L}. In Section~\ref{scomonoid} we give a
self-contained presentation of this comonoid.

In the following section we then show that $\E$ acts on the comonoid,
and hence defines a natural transformation from the forgetful functor
to itself. Considered as an element of a completion of $U_q\g$ this
transformation is a $1$-cochain with coboundary $\E$. Pushing the
analysis further we then conclude that $\E=1$.

In Section~\ref{SecTwist} we reformulate the main result as a statement
about uniqueness of the Drinfeld twist, and in Section~\ref{sdirac} we
apply this to show that the quantum Dirac operator is uniquely defined
up to unitary equivalence.

We end the paper with two appendices. In Appendix A we prove the
essentially known result that any group-like element affiliated with
$W^*(G)$ belongs to the complexification of $G$. This is used in the
main text to show that the group of central group-like elements of the
completion of~$U_q\g$ is isomorphic to the center of $G$. In Appendix B
we provide a short proof of our main result in the formal deformation
setting following Drinfeld's arguments for $3$-cocycles.


\bigskip

\section{Cohomology of quantum groups}
\label{seccohom}

Let $(A,\Delta)$ be a discrete bialgebra in the sense of~\cite{NT3}.
Therefore $A\cong \oplus_{\lambda\in\Lambda}\End(V_\lambda)$ as an
algebra, and
$$
\Delta\colon A\to M(A\otimes A)\cong\prod_{\lambda,\mu}\End(V_\lambda\otimes V_\mu)
$$
is a nondegenerate homomorphism satisfying coassociativity and
which comes with a counit $\eps$.

Adapting the usual definition of cohomology for Hopf algebras, see e.g.
Section 2.3 in \cite{Ma}, define an operator
$$
\partial\colon M(A^{\otimes n})^\times\to M(A^{\otimes(n+1)})^\times,
$$
where the superindex $\times$ denotes invertible elements, by
$$\partial\chi=(\Delta_0 (\chi)\Delta_2
(\chi)\dots )(\Delta_1 (\chi^{-1})\Delta_3 (\chi^{-1})\dots),$$
where $\Delta_i
=\iota\otimes\dots\otimes\iota\otimes\Delta\otimes\iota\otimes\dots\otimes\iota$
with $\Delta$ in the $i$th position for $0<i<n+1$, and $\Delta_0(\chi)
=1\otimes\chi$, $\Delta_{n+1}(\chi)=\chi\otimes 1$. In particular, if
$u\in M(A)^\times$ and $\E\in M(A\otimes A)^\times$, then we have
$$\partial u=(u\otimes u)\Delta(u)^{-1},\ \ \ \partial\E =(1\otimes\E)
(\iota\otimes\Delta)(\E)(\Delta\otimes\iota)(\E^{-1})(\E^{-1}\otimes
1).$$ An $n$-cochain $\chi\in M(A^{\otimes n})^\times$ is a called a
cocycle if $\partial\chi =1$, and it is called a coboundary if $\chi$
belongs to the image of $\partial$. In general, an $n$-coboundary is not
necessarily a cocycle, but this is the case for $n=2$.

Two $2$-cochains $\E,\F$ are said to be cohomologous if there exists a
$1$-cochain $u$ such that
$$\E=(u\otimes u)\F\Delta(u)^{-1}.$$
Then $\partial\E =(u\otimes u\otimes u)\partial\F(u^{-1}\otimes
u^{-1}\otimes u^{-1})$, so if $\F$ is a cocycle, then so is $\E$, and
we also conclude that $\partial^2 u=1$, although in general
$\partial^2\neq 1$. The set of cohomology classes of $2$-cocycles is
denoted by~$H^2 (A;\C^*)$; this is just a set, the product of two
$2$-cocycles is not necessarily a cocycle.

\medskip

We say that a $2$-cocycle $\E\in
M(A\otimes A)$ is
\begin{itemize}
\item {\bf invariant}, if $[\E,\Delta(a)]=0$ for all $a\in A$;
    \item {\bf symmetric}, if $A$ is quasitriangular with
        $R$-matrix $\RR$ and $\RR\E=\E_{21}\RR$.
\end{itemize}

If a $2$-cocycle $\E$ is invariant and symmetric, then the cohomologous cocycle
$\F=(v\otimes v)\E\Delta(v)^{-1}$ is also invariant and symmetric whenever $v$ is central.

\medskip

If $A$ is a discrete $*$-bialgebra, so that $A$ completes to a
C$^*$-algebra and $\Delta$ is a $*$-homomorphism, it makes more sense
to consider only unitary cochains and define $H^2 (A;\T)$. The following simple lemma shows that the canonical map $H^2 (A;\T)\to H^2 (A;\C^*)$ is injective.

\begin{lemma}
\label{lempolar} Suppose $(A, \Delta)$ is a discrete $*$-bialgebra.
Consider two unitary $2$-cochains $\E,\F$ such that $\E =(u\otimes
u)\F\Delta(u)^{-1}$ for a $1$-cochain $u$. Then $\E =(v\otimes
v)\F\Delta(v)^{-1}$, where $v$ is the unitary part in the polar
decomposition $u=v|u|$.
\end{lemma}

\bp It is sufficient to show that $(|u|\otimes|u|)\F=\F\Delta(|u|)$, or
since $|u|=\sqrt{u^*u}$, that $$(u^*u\otimes u^*u)\F=\F\Delta(u^*u),$$
and this is immediate from
$$1=\E^*\E=\Delta(u^{-1})^*\F^*(u^*\otimes u^*)(u\otimes u)\F\Delta(u^{-1}).$$
\ep

For invariant $2$-cocycles it makes also sense to consider the polar
decomposition.

\begin{lemma}
Suppose $(A, \Delta)$ is a (quasitriangular) discrete $*$-bialgebra.
Let $\E$ be a (symmetric) invariant $2$-cocycle, and $\E=\F|\E|$ be the
polar decomposition. Then both $\F$ and $|\E|$ are (symmetric)
invariant $2$-cocycles.
\end{lemma}

\bp A somewhat more general statement is proved in
\cite[Proposition~1.3]{NT3}. Briefly, observe that (symmetric)
invariant cocycles form a group which is in addition closed under
involution (recall that for quasitriangular $*$-bialgebras we require
$\RR^*=\RR_{21}$). It follows that $\E^*\E$ is again such a cocycle.
Taking the square root, we conclude that $|\E|$ and $\F=\E|\E|^{-1}$
are (symmetric) invariant cocycles as well. \ep

We end this section with a categorical perspective which is convenient
to keep in mind. Consider the category $A\Mod$ of finite dimensional
nondegenerate $A$-modules, and let $F\colon A\Mod\to\Vect$ be the
forgetful functor. Then $M(A)$ as an algebra is identified with the
algebra $\Nat(F)$ of natural transformations from $F$ to itself.

An invertible element $\E\in M(A\otimes A)$ defines a natural
isomorphism
$$
F_2\colon F(U)\otimes F(V)\xrightarrow{\E^{-1}} F(U\otimes V).
$$
Assume that $(\eps\otimes\iota)(\E)=(\iota\otimes\eps)(\E)=1$.
Then $\E$ is a cocycle if $(F,F_2,F_0=\iota)$ is a tensor functor, that
is, the diagram
$$
\xymatrix{
F(U)\otimes F(V)\otimes F(W) \ar[rr]^{\E^{-1}\otimes1}\ar[d]_{1\otimes\E^{-1}}&
& F(U\otimes V)\otimes F(W) \ar[d]^{(\Delta\otimes\iota)(\E^{-1})}\\
F(U)\otimes F(V\otimes W) \ar[rr]_{(\iota\otimes\Delta)(\E^{-1})} & &
F(U\otimes V\otimes W)}
$$
commutes. (We remark that a $2$-cocycle such that
$(\eps\otimes\iota)(\E)=(\iota\otimes\eps)(\E)=1$ is called
counital. Any $2$-cocycle is cohomologous to a counital one, since by
applying $\eps$ to the middle term of the cocycle identity we
conclude that
$(\eps\otimes\iota)(\E)=(\iota\otimes\eps)(\E)=(\eps\otimes\eps)(\E)1$.
If $\E$ is not counital, to define a tensor functor we just have to put
$F_0=(\eps\otimes\eps)(\E)$ instead of $F_0=\iota$.)

Two $2$-cocycles are cohomologous if the corresponding tensor functors $A\Mod\to\Vect$
are naturally isomorphic.

\smallskip

A cocycle $\E$ is invariant if it
defines a natural morphism $U\otimes V\to U\otimes V$ in $A\Mod$. In this case the identity functor $F\colon A\Mod\to A\Mod$
becomes a tensor functor with $F_2=\E^{-1}$. The invariant cocycle $\E$ is symmetric if and only if this tensor functor is braided, that is, the diagram
$$
\xymatrix{
F(U)\otimes F(V) \ar[r]^{\E^{-1}}\ar[d]_{\sigma}& F(U\otimes V) \ar[d]^{F(\sigma)}\\
F(V)\otimes F(U) \ar[r]_{\E^{-1}} & F(V\otimes U)}
$$
commutes, where $\sigma=\Sigma\RR$ is the braiding. For two invariant cocycles~$\E$ and~$\F$, there exists a central element $c$ such that $\E=(c\otimes c)\F\Delta(c)^{-1}$ if and only if the corresponding tensor functors $A\Mod\to A\Mod$ are naturally isomorphic.

\smallskip

For a discussion of the Drinfeld associator, which is a $3$-cocycle, see
Section \ref{SecTwist}.

\bigskip

\section{Main result}\label{smain}

Let $G$ be a simply connected semisimple compact Lie group, $\g$ its
complexified Lie algebra. Denote by~$\widehat{\C[G]}$ the discrete
bialgebra of matrix coefficients of finite dimensional representations
of~$G$ with convolution product. It is a quasitriangular discrete
$*$-bialgebra with $R$-matrix $\RR=1$. We write~$\U(G)$ instead of
$M(\widehat{\C[G]})$. It is the algebra of closed densely defined
operators affiliated with the von Neumann algebra $W^*(G)$ of $G$, and
it contains the universal enveloping algebra $U\g$. We denote
by~$\Dhat$ and $\hat\eps$ the comultiplication and the counit on
$\U(G)$.

\smallskip

Let $\h\subset\g$ be the Cartan subalgebra defined by a maximal torus
in $G$. Fix a system $\{\alpha_1,\dots,\alpha_r\}$ of simple roots. Let
$\omega_1,\dots,\omega_r$ be the fundamental weights. The weight and
root lattices are denoted by $P$ and $Q$, respectively. Let
$(a_{ij})_{1\le i,j\le r}$ be the Cartan matrix and $d_1,\dots,d_r$ be
the coprime positive integers such that $(d_ia_{ij})_{i,j}$ is
symmetric. Define a bilinear form on $\h^*$ by $(\alpha_i,\alpha_j)
=d_ia_{ij}$. Let $h_i\in\h$ be such that $\alpha_j(h_i)=a_{ij}$. For
$\lambda\in P$ we shall write $\lambda(i)$ instead of $\lambda(h_i)$.
Therefore $\lambda(1),\dots,\lambda(r)$ are the coefficients of
$\lambda$ in the basis $\omega_1,\dots,\omega_r$.

\smallskip

For $q\in\C^*$ not a root of unity consider the quantized universal
enveloping algebra $U_q\g$ generated by elements $E_i$, $F_i$, $K_i$,
$K_i^{-1}$, $1\le i\le r$, satisfying the relations
$$
K_iK_i^{-1}=K_i^{-1}K_i=1,\ \ K_iK_j=K_jK_i,\ \
K_iE_jK_i^{-1}=q_i^{a_{ij}}E_j,\ \
K_iF_jK_i^{-1}=q_i^{-a_{ij}}F_j,
$$
$$
E_iF_j-F_jE_i=\delta_{ij}\frac{K_i-K_i^{-1}}{q_i-q_i^{-1}},
$$
$$
\sum^{1-a_{ij}}_{k=0}(-1)^k\begin{bmatrix}1-a_{ij}\\
k\end{bmatrix}_{q_i} E^k_iE_jE^{1-a_{ij}-k}_i=0,\ \
\sum^{1-a_{ij}}_{k=0}(-1)^k\begin{bmatrix}1-a_{ij}\\
k\end{bmatrix}_{q_i} F^k_iF_jF^{1-a_{ij}-k}_i=0,
$$
where $\displaystyle\begin{bmatrix}m\\
k\end{bmatrix}_{q_i}=\frac{[m]_{q_i}!}{[k]_{q_i}![m-k]_{q_i}!}$,
$[m]_{q_i}!=[m]_{q_i}[m-1]_{q_i}\dots [1]_{q_i}$,
$\displaystyle[n]_{q_i}=\frac{q_i^n-q_i^{-n}}{q_i-q_i^{-1}}$ and
$q_i=q^{d_i}$. This is a Hopf algebra with coproduct $\Dhat_q$ and
counit $\hat\eps_q$ defined by
$$
\Dhat_q(K_i)=K_i\otimes K_i,\ \
\Dhat_q(E_i)=E_i\otimes1+ K_i\otimes E_i,\ \
\Dhat_q(F_i)=F_i\otimes K_i^{-1}+1\otimes F_i,
$$
$$
\hat\eps_q(E_i)=\hat\eps_q(F_i)=0,\ \ \hat\eps_q(K_i)=1.
$$

If $V$ is a finite dimensional $U_q\g$-module and $\lambda\in P$,
denote by $V(\lambda)$ the space of vectors $v\in V$ of weight
$\lambda$, so that $K_iv=q_i^{\lambda(i)}v$ for all $i$. Recall that
$V$ is called admissible if $V=\oplus_{\lambda\in P}V(\lambda)$.
Consider the tensor category of finite dimensional admissible
$U_q\g$-modules. It is a semisimple category with simple objects
indexed by dominant integral weights $\lambda\in P_+$. For each
$\lambda\in P_+$ we fix an irreducible $U_q\g$-module $V_\lambda$ with
highest weight $\lambda$. Denote by $\widehat{\C[G_q]}$ the discrete
bialgebra defined by our category, so
$\widehat{\C[G_q]}\cong\oplus_{\lambda\in P_+}\End(V_\lambda)$. Denote
by $\U(G_q)$ the multiplier algebra $M(\widehat{\C[G_q]})$. We shall
write $\U(G_q\times G_q)$ instead of $M(\widehat{\C[G_q]}\otimes
\widehat{\C[G_q]})$.

The discrete bialgebra $\widehat{\C[G_q]}$ is quasitriangular. The
$R$-matrix $\RR_\hbar$ depends on the choice of $\hbar\in\C$ such that
$q=e^{\pi i\hbar}$. From now on we will write $q^x$ instead of $e^{\pi
i\hbar x}$, provided the choice of~$\hbar$ is clear from the context.
The $R$-matrix $\RR_\hbar$ can can be defined by an explicit formula,
see e.g. \cite[Theorem~8.3.9]{CP}, but for us it will be enough to
remember that it is characterized by the following two properties:
\begin{itemize}
\item $\Dhat^{op}_q=\RR_\hbar\Dhat_q(\cdot)\RR^{-1}_\hbar$;
\item if $U$ is a module with a lowest weight vector $\zeta$ of
    weight $\lambda$, so $\zeta\in U(\lambda)$ and $F_i\zeta=0$ for
    all $i$, and $V$ is a module with a highest weight vector $\xi$
    of weight $\mu$, so $\xi\in V(\mu)$ and $E_i\xi=0$ for all~$i$,
    then
\begin{equation} \label{eRmat}
\RR_\hbar(\zeta\otimes\xi)
=q^{(\lambda,\mu)}\zeta\otimes\xi.
\end{equation}
\end{itemize}
This indeed characterizes $\RR_\hbar$, since if $U$ and $V$ are
irreducible then $\zeta\otimes\xi$ is a cyclic vector in $U\otimes V$.

We denote by $\CC(\g,\hbar)$ the braided monoidal category of
admissible finite dimensional $U_q\g$-modules with braiding
$\sigma=\Sigma\RR_\hbar$.

\medskip

We are now ready to formulate our main result.

\begin{theorem} \label{thmMain}
Assume $q\in\C^*$ is not a nontrivial root of unity. Then any symmetric
invariant $2$-cocycle $\E\in\U(G_q\times G_q)$ is the coboundary of a
central element in $\U(G_q)$.
\end{theorem}

We will assume $q\ne1$, leaving the case $q=1$, which requires minor,
mostly notational, modifications, to the reader. In fact, the proof of
\cite[Theorem~11]{Wa1} and \cite[Lemma~29(1)]{Wa1} show that for $q=1$
the result is true for any connected compact group $G$, at least for
unitary cocycles.

\smallskip

For each $\mu\in P_+$ fix a highest weight vector $\xi_\mu\in V_\mu$.
We identify $V_0$ with $\C$ so that $\xi_0=1$.

For $\mu,\eta\in P_+$, define a morphism
$$
T_{\mu,\eta}\colon V_{\mu+\eta}\to V_\mu\otimes V_\eta \ \ \hbox{by}\ \
\xi_{\mu+\eta}\mapsto\xi_\mu\otimes \xi_\eta.
$$
The image of $T_{\mu,\eta}$ is the irreducible isotypic component of
$V_\mu\otimes V_\eta$ with highest weight $\mu+\eta$. Since~$\E$ is
invariant, the action of $\E$ on this image is by a scalar, so there
exists $\E(\mu,\eta)\in\C^*$ such that
$$
\E T_{\mu,\eta}=\E(\mu,\eta)T_{\mu,\eta}.
$$

\begin{lemma}
The map $P_+\times P_+\to \C^*$, $(\mu,\eta)\mapsto\E(\mu,\eta)$, is a
symmetric $2$-cocycle, that is,
$$
\E(\mu,\eta)=\E(\eta,\mu)\ \ \text{and}\ \
\E(\mu,\eta)\E(\mu+\eta,\nu)
=\E(\eta,\nu)\E(\mu,\eta+\nu).
$$
\end{lemma}

\bp That it is a cocycle follows from the identity
\begin{equation*}\label{eTT}
(T_{\mu,\eta}\otimes\iota)T_{\mu+\eta,\nu}
=(\iota\otimes T_{\eta,\nu})T_{\mu,\eta+\nu}
\end{equation*}
by applying the operator $(\E\otimes1)(\Dhat_q\otimes\iota)(\E)$ to the
left side and the same operator $(1\otimes\E)(\iota\otimes\Dhat_q)(\E)$
to the right side.

To see that the cocycle is symmetric, notice that the braiding $\sigma$
maps the image of $T_{\mu,\eta}$, which is the irreducible isotypic
component with highest weight $\mu+\eta$, onto the image of
$T_{\eta,\mu}$ (in fact, one can show that $\sigma
T_{\mu,\eta}=q^{(\mu,\eta)}T_{\eta,\mu}$). Since $\sigma\E=\E\sigma$
this gives the result. \ep

It is well-known that any symmetric $2$-cocycle on $P_+$ is a
coboundary, see e.g. \cite[Lemma~4.2]{NT3}, that is, there exist
$c(\mu)\in\C^*$ such that
$$
\E(\mu,\eta)=c(\mu+\eta)c(\mu)^{-1}c(\eta)^{-1}.
$$
The numbers $c(\mu)$, $\mu\in P_+$, define an invertible element $c$ in
the center of $\U(G_q)$. Then replacing~$\E$ by $(c\otimes
c)\E\Dhat_q(c)^{-1}$ we get a new symmetric invariant $2$-cocycle which
is cohomologous to $\E$ via a central element and is such that the
corresponding $2$-cocycle on $P_+$ is trivial. In other words, without
loss of generality we may assume that
\begin{equation} \label{eHigh}
\E(\mu,\eta)=1\ \ \hbox{for all}\ \ \mu,\eta\in P_+.
\end{equation}
This in particular implies that $\E$ is counital, since
$(\hat\eps_q\otimes\iota)(\E)$ acts on $V_\mu$ as multiplication by
$\E(0,\mu)$.

\medskip

Fix $i$, $1\le i\le r$. Let $\mu,\eta\in P_+$ be such that
$\mu(i),\eta(i)\ge1$. Then the space $(V_\mu\otimes
V_\eta)(\mu+\eta-\alpha_i)$ is $2$-dimensional, spanned by
$F_i\xi_\mu\otimes\xi_\eta$ and $\xi_\mu\otimes F_i\eta$. This space
has a unique, up to a scalar, vector killed by $E_i$, namely,
$$
[\mu(i)]_{q_i}\xi_\mu\otimes F_i\xi_\eta
-q_i^{\mu(i)}[\eta(i)]_{q_i}F_i\xi_\mu\otimes\xi_\eta.
$$
In other words, the isotypic component of $V_\mu\otimes V_\eta$ with
highest weight $\mu+\eta-\alpha_i$ is the image of the morphism
\begin{equation*} \label{etau}
\tau_{i;\mu,\eta}\colon V_{\mu+\eta-\alpha_i}\to V_\mu\otimes
V_\eta, \ \ \xi_{\mu+\eta-\alpha_i}\mapsto
[\mu(i)]_{q_i}\xi_\mu\otimes F_i\xi_\eta
-q_i^{\mu(i)}[\eta(i)]_{q_i}F_i\xi_\mu\otimes\xi_\eta.
\end{equation*}
The action of $\E$ on this image is by a scalar, so there exists
$\E_i(\mu,\eta)\in\C^*$ such that
$$
\E \tau_{i;\mu,\eta}=\E_i(\mu,\eta)\tau_{i;\mu,\eta}.
$$

\begin{lemma}
Assume the cocycle $\E$ satisfies condition \eqref{eHigh}. Then, for
fixed $i$, the numbers $\E_i(\mu,\eta)$ do not depend on $\mu$ and
$\eta$ with $\mu(i),\eta(i)\ge1$.
\end{lemma}

\bp Consider the module $V_\mu\otimes V_\eta\otimes V_\nu$. The
isotypic component corresponding to $\mu+\eta+\nu-\alpha_i$ has
multiplicity two, and is spanned by the images of $(\iota\otimes
T_{\eta,\nu})\tau_{i;\mu,\eta+\nu}$ and
$(\iota\otimes\tau_{i;\eta,\nu})T_{\mu,\eta+\nu-\alpha_i}$, as well as
by the images of $(T_{\mu,\eta}\otimes\iota)\tau_{i;\mu+\eta,\nu}$ and
$(\tau_{i;\mu,\eta}\otimes\iota)T_{\mu+\eta-\alpha_i,\nu}$. These maps
are related by the following identities:
\begin{equation}\label{etauT1}
[\eta(i)]_{q_i}(T_{\mu,\eta}\otimes\iota)\tau_{i;\mu+\eta,\nu}
-[\nu(i)]_{q_i}(\tau_{i;\mu,\eta}\otimes\iota)
T_{\mu+\eta-\alpha_i,\nu}=[\mu(i)+\eta(i)]_{q_i}
(\iota\otimes\tau_{i;\eta,\nu})
T_{\mu,\eta+\nu-\alpha_i},
\end{equation}
\begin{equation}\label{etauT2}
[\eta(i)+\nu(i)]_{q_i}(\tau_{i;\mu,\eta}\otimes\iota)
T_{\mu+\eta-\alpha_i,\nu}=[\eta(i)]_{q_i}
(\iota\otimes T_{\eta,\nu})\tau_{i;\mu,\eta+\nu}-
[\mu(i)]_{q_i}(\iota\otimes\tau_{i;\eta,\nu})
T_{\mu,\eta+\nu-\alpha_i}.
\end{equation}
These identities are checked by applying both sides to the highest
weight vector $\xi_{\mu+\eta+\nu-\alpha_i}$, and using that the $T$'s are
module maps, so that for example
$$
T_{\mu,\eta}F_i=(F_i\otimes K_i^{-1})T_{\mu,\eta}+(1\otimes F_i)T_{\mu,\eta}.
$$

Applying
$(\E\otimes1)(\Dhat_q\otimes\iota)(\E)=(1\otimes\E)(\iota\otimes\Dhat_q)(\E)$
to \eqref{etauT1} and using that $\E T=T$ by \eqref{eHigh}, we get
\begin{equation*}
\begin{split}
\E_i(\mu+\eta,\nu)[\eta(i)]_{q_i}(T_{\mu,\eta}\otimes\iota)\tau_{i;\mu+\eta,\nu}
-\E_i(\mu,\eta)[\nu(i)]_{q_i}(\tau_{i;\mu,\eta}\otimes\iota)
T_{\mu+\eta-\alpha_i,\nu}\\=\E_i(\eta,\nu)[\mu(i)+\eta(i)]_{q_i}
(\iota\otimes\tau_{i;\eta,\nu})
T_{\mu,\eta+\nu-\alpha_i}
\end{split}
\end{equation*}
Since $(T_{\mu,\eta}\otimes\iota)\tau_{i;\mu+\eta,\nu}$ and
$(\tau_{i;\mu,\eta}\otimes\iota) T_{\mu+\eta-\alpha_i,\nu}$ are
linearly independent, together with \eqref{etauT1} this implies that
$$
\E_i(\mu+\eta,\nu)=\E_i(\mu,\eta)=\E_i(\eta,\nu).
$$
Now for arbitrary $\mu,\eta,\tilde\mu,\tilde\eta$, applying the last
identity twice, we get
$\E_i(\mu,\eta)=\E_i(\eta,\tilde\mu)=\E_i(\tilde\mu,\tilde\eta). $\ep

Define a homomorphism $\chi\colon Q\to\C^*$ by letting
$\chi(\alpha_i)=\E_i(\mu,\eta)^{-1}$ for $\mu,\eta\in P_+$ with
$\mu(i),\eta(i)\ge1$, $1\le i\le r$. Extend $\chi$ to a homomorphism
$P\to\C^*$. The restriction of $\chi$ to $P_+$ defines a central
element $c$ of~$\U(G_q)$ such that
$$
(c\otimes c)\Dhat_q(c)^{-1}\tau_{i;\mu,\eta}
=\chi(\mu)\chi(\eta)\chi(\mu+\eta-\alpha_i)^{-1}\tau_{i;\mu,\eta}
=\chi(\alpha_i)\tau_{i;\mu,\eta}=\E_i(\mu,\eta)^{-1}\tau_{i;\mu,\eta}.
$$
Thus replacing $\E$ by the cohomologous cocycle $(c\otimes
c)\E\Dhat_q(c)^{-1}$ we get a symmetric invariant $2$-cocycle, which
we again denote by $\E$, such that
\begin{equation}\label{eTauNorm}
\E_i(\mu,\eta)=1\ \ \hbox{for all}\ \ 1\le i\le r\ \ \hbox{and}\ \ \mu,\eta\in P_+\ \
\hbox{with}\ \ \mu(i),\eta(i)\ge1.
\end{equation}
Note that condition \eqref{eHigh} for this new cocycle is still
satisfied, since $\chi$ is a homomorphism on $P_+$.

\medskip

From now on we can and will assume that the symmetric invariant
$2$-cocycle $\E$ satisfies properties~\eqref{eHigh} and
\eqref{eTauNorm}. We will see later that this already implies that
$\E=1$. But to show this we have to make a rather long detour and first
prove that $\E$ is the coboundary of a central element. In the
remaining part of the section we will show that for $G=\SU(2)$ this can
be avoided.

\medskip

Recall that for $G=\SU(2)$ the weight lattice $P$ is identified with
$\frac{1}{2}\Z$ and the root lattice with $\Z$. For
$s\in\frac{1}{2}\N$, we have $V_{1/2}\otimes V_s \cong V_{s+1/2}\oplus
V_{s-1/2}$. Therefore conditions~\eqref{eHigh} and \eqref{eTauNorm}
imply that $\E$ acts trivially on $V_{1/2}\otimes V_s$.

Now for $s,t\ge 1/2$ consider the morphism $T_{1/2,s}\otimes\iota\colon
V_{s+1/2}\otimes V_t\to V_{1/2}\otimes V_s\otimes V_t$ and compute
\begin{align*}
(T_{1/2,s}\otimes\iota)\E
&=(\Dhat_q\otimes\iota)(\E)(T_{1/2,s}\otimes\iota)\\
&=(1\otimes\E)(\iota\otimes\Dhat_q)(\E)
(\E^{-1}\otimes1)(T_{1/2,s}\otimes\iota)\\
&=(1\otimes\E)
(T_{1/2,s}\otimes\iota),
\end{align*}
since $\E$ acts trivially on $V_{1/2}\otimes V$ for any $V$. It follows
that if $\E$ acts trivially on $V_s\otimes V_t$, it acts trivially on
$V_{s+1/2}\otimes V_t$. Therefore an induction argument shows that $\E$
acts trivially on $V_s\otimes V_t$ for all $s$ and $t$, so $\E=1$.

\smallskip

For general $G$ one can similarly show that it suffices to check that
$\E$ acts trivially on $V_{\omega_i}\otimes V_\mu$, but we don't know
whether it is possible to check the latter property directly using
conditions \eqref{eHigh} and \eqref{eTauNorm}.

\bigskip

\section{Comonoid representing the canonical fiber functor}
\label{scomonoid}

Consider the automorphism $\theta$ of $U_q\g$ defined by
$$
\theta(E_i)=F_i,\ \ \theta(F_i)=E_i,\ \ \theta(K_i)=K_i^{-1}.
$$
Observe that $\Dhat_q\theta=(\theta\otimes\theta)\Dhat_q^{op}$.
For every $U_q\g$-module $V$ define a module $\bar V$ which coincides
with~$V$ as a vector space, but the action of $U_q\g$ is given by
$$
X\bar v=\overline{\theta(X)v},
$$
where $\bar v$ means the vector $v\in V$ considered as an element of
$\bar V$. Notice that $\bar\xi_\mu$ is a lowest weight vector of
weight~$-\mu$.

Denote by $\bar\mu$ the weight $-w_0\mu$, where $w_0$ is the longest
element in the Weyl group. The involution $\lambda\mapsto\bar\lambda$
on $P$ defines an involution on the index set $\{1,\dots,r\}$, so that
$\bar\alpha_i=\alpha_{\bar i}$ and $\bar\omega_i=\omega_{\bar i}$. It
is known that the lowest weight of $V_\mu$ is $-\bar\mu$. It follows
that $\bar V_\mu\cong V_{\bar\mu}$.

For each $\mu\in P_+$ there exists a unique up to a scalar morphism
$\bar V_\mu\otimes V_\mu\to V_0=\C$. Namely, define
$$
S_\mu\colon \bar V_\mu\otimes V_\mu\to\C,\ \
\bar\xi_\mu\otimes\xi_\mu\mapsto1,
$$
see e.g. \cite[Proposition 25.1.4]{L}.

For $\mu,\eta\in P_+$, define a morphism
$$
\bar T_{\mu,\eta}\colon \bar V_{\mu+\eta}\to \bar V_\mu\otimes \bar V_\eta \ \
\hbox{by}\ \
\bar\xi_{\mu+\eta}\mapsto\bar\xi_\mu\otimes \bar\xi_\eta.
$$

For $\lambda\in P$ and $\mu,\eta\in P_+$ such that $\lambda+\mu\in P_+$
consider the morphism
$$
\tr^{\eta}_{\mu,\lambda+\mu}\colon \bar V_{\mu+\eta}\otimes
V_{\lambda+\mu+\eta}\to \bar V_{\mu}\otimes V_{\lambda+\mu},\ \
\bar\xi_{\mu+\eta}\otimes\xi_{\lambda+\mu+\eta}
\mapsto\bar\xi_{\mu}\otimes\xi_{\lambda+\mu}.
$$
Since $\bar\xi_{\mu+\eta}\otimes\xi_{\lambda+\mu+\eta}$ is a cyclic
vector, its image completely determines the morphism, if it exists. To
show existence, rewrite this morphism as the composition
\begin{equation*} \label{etr0}
\bar V_{\mu+\eta}\otimes
V_{\lambda+\mu+\eta} \xrightarrow{\bar T_{\mu,\eta}\otimes
T_{\eta,\lambda+\mu}} \bar V_{\mu}\otimes \bar V_{\eta}\otimes
V_\eta\otimes V_{\lambda+\mu} \xrightarrow{\iota\otimes
S_\eta\otimes\iota}\bar V_{\mu}\otimes V_{\lambda+\mu}.
\end{equation*}

Using the morphisms $\tr$ define the inverse limit $U_q\g$-module
$$
M_\lambda=\lim_{\xleftarrow[\mu]{}}\bar V_\mu\otimes
V_{\lambda+\mu}.
$$
We consider $M_\lambda$ as a topological $U_q\g$-module with a base of
neighborhoods of zero formed by the kernels of the canonical morphisms
$M_\lambda\to \bar V_{\mu}\otimes V_{\lambda+\mu}$. Observe that
$\tr^{\eta}_{\mu,\lambda+\mu}$ is surjective since its image contains
the cyclic vector $\bar\xi_{\mu}\otimes\xi_{\lambda+\mu}\in \bar
V_\mu\otimes V_{\lambda+\mu}$. It follows that the morphisms
$M_\lambda\to \bar V_{\mu}\otimes V_{\lambda+\mu}$ are surjective.
Hence, if $V$ is a $U_q\g$-module with discrete topology, then any
continuous morphism $M_\lambda\to V$ factors through $\bar
V_{\mu}\otimes V_{\lambda+\mu}$ for some~$\mu$, so that the space
$\Hom_{U_q\g}(M_\lambda,V)$ of such morphisms is the inductive limit of
$\Hom_{U_q\g}(\bar V_{\mu}\otimes V_{\lambda+\mu},V)$.

Recall, see e.g. \cite[Proposition 23.3.10]{L}, that if $V$ is an
admissible finite dimensional $U_q\g$-module and $\lambda$ an integral
weight then the map
\begin{equation} \label{emult}
\Hom_{U_q\g}(\bar V_{\mu}\otimes V_{\lambda+\mu},V)\to V(\lambda),\ \
f\mapsto f(\bar\xi_{\mu}\otimes\xi_{\lambda+\mu}),
\end{equation}
is an isomorphism for sufficiently large dominant integral weights
$\mu$. 
In particular, for any $V\in\CC(\g,\hbar)$ the maps \eqref{emult}
induce a linear isomorphism
$$
\Hom_{U_q\g}(M_\lambda,V)\to V(\lambda).
$$
Therefore the topological $U_q\g$-module $M=\oplus_{\lambda\in P}
M_\lambda$ represents the forgetful functor $\CC(\g,\hbar)\to\Vect$.
Let
$$
\eta_V\colon \Hom_{U_q\g}(M,V)\to V
$$
be the canonical isomorphism.

Our next goal is to define a comonoid structure on $M$. Define
$$
M_{\lambda_1}\hat\otimes M_{\lambda_2}
=\lim_{\xleftarrow[\mu_1,\mu_2]{}}(\bar V_{\mu_1}\otimes
V_{\lambda_1+\mu_1})\otimes (\bar V_{\mu_2}\otimes V_{\lambda_2+\mu_2})
$$
and then
$$
M\hat\otimes M=\prod_{\lambda_1,\lambda_2\in P}
M_{\lambda_1}\hat\otimes M_{\lambda_2}.
$$
Higher tensor powers of $M$ are defined similarly. We want to define a
morphism
$$
\delta\colon M\to M\hat\otimes M.
$$
The restriction of $\delta$ to $M_\lambda$ composed with the projection
$M\hat\otimes M\to M_{\lambda_1}\hat\otimes M_{\lambda_2}$ will be
nonzero only if $\lambda=\lambda_1+\lambda_2$, so $\delta$ is
determined by maps
$$
\delta_{\lambda_1,\lambda_2}\colon
M_{\lambda_1+\lambda_2} \to M_{\lambda_1}\hat\otimes
M_{\lambda_2}.
$$
We define these morphisms using the morphisms
$$
m_{\mu,\eta,\lambda_1,\lambda_2}\colon
\bar V_{\mu+\eta}\otimes V_{\lambda_1+\lambda_2+\mu+\eta}\to
\bar V_{\mu}\otimes  V_{\lambda_1+\mu}\otimes
\bar V_{\eta}\otimes V_{\lambda_2+\eta}
$$
mapping $\bar\xi_{\mu+\eta}\otimes \xi_{\lambda_1+\lambda_2+\mu+\eta}$
onto $\bar\xi_{\mu}\otimes\xi_{\lambda_1
+\mu}\otimes\bar\xi_{\eta}\otimes\xi_{\lambda_2+\eta}$. Since
$\bar\xi_{\mu+\eta}\otimes \xi_{\lambda_1+\lambda_2+\mu+\eta}$ is a
cyclic vector, such a morphism is unique if it exists, and to show its
existence we rewrite it, using property \eqref{eRmat} of the
$R$-matrix, as the composition
\begin{equation*} \label{em}
\begin{split}
\bar V_{\mu+\eta}\otimes V_{\lambda_1+\lambda_2+\mu+\eta}
&\xrightarrow{\bar T_{\mu,\eta}\otimes
T_{\lambda_1+\mu,\lambda_2+\eta}}\bar V_{\mu}\otimes
\bar V_{\eta}\otimes V_{\lambda_1+\mu}
\otimes V_{\lambda_2+\eta}\\
&\xrightarrow{q^{(\lambda_1+\mu,\eta)}(\iota\otimes \sigma
\otimes\iota)}\bar V_{\mu}\otimes  V_{\lambda_1+\mu}\otimes
\bar V_{\eta}\otimes V_{\lambda_2+\eta}.
\end{split}
\end{equation*}

The morphisms $m$ are consistent with the morphisms $\tr$ defining the
inverse limits, that is,
$$
(\tr^{\nu}_{\mu,\lambda_1+\mu}\otimes
\tr^{\omega}_{\eta,\lambda_2+\eta})
m_{\mu+\nu,\eta+\omega,\lambda_1,\lambda_2}
=m_{\mu,\eta,\lambda_1,\lambda_2}
\tr^{\nu+\omega}_{\mu+\eta,\lambda_1+\lambda_2+\mu+\eta}.
$$
Hence they define morphisms $\delta_{\lambda_1,\lambda_2}\colon
M_{\lambda_1+\lambda_2} \to M_{\lambda_1}\hat\otimes M_{\lambda_2}$.

Using the morphisms $\delta_{\lambda_1,\lambda_2}$ we can in an obvious
way define morphisms
$$
(\delta\otimes\iota)\delta,
(\iota\otimes\delta)\delta\colon M\to
M\hat\otimes M\hat\otimes M.
$$
We also introduce a morphism $\eps\colon M\to\C$ by requiring it to be
nonzero only on $M_0$, where we set it to be the canonical morphism
$M_0\to \bar V_0\otimes V_0=\C$, so that $\eps\colon M_0\to\C$ is
determined by the morphisms
$$
\tr^{\mu}_{0,0}=S_\mu\colon \bar V_{\mu}\otimes V_\mu\to\C.
$$

\begin{proposition} \label{pComonoid0}
The triple $(M,\delta,\eps)$ is a comonoid representing the canonical
fiber functor $\CC(\g,\hbar)\to \Vect$, that is,
$$(\delta\otimes\iota)\delta
=(\iota\otimes\delta)\delta,\ \ (\eps\otimes\iota)\delta=\iota
=(\iota\otimes\eps)\delta,
$$
and for all $U,V\in\CC(\g,\hbar)$ the following diagram commutes:
$$
\xymatrix{
\Hom_{U_q\g}(M,U)\otimes \Hom_{U_q\g}(M,V)
\ar[rr]^{\mbox{\quad\quad\quad\quad\quad}\eta_U\otimes\eta_V}\ar[d]&
& U\otimes V\ar@{=}[d]\\
\Hom_{U_q\g}(M,U\otimes V) \ar[rr]^{\mbox{\quad\quad\quad\quad\quad}
\eta_{U\otimes V}} &  &
U\otimes V},
$$
where the left vertical arrow is given by $f\otimes g\mapsto (f\otimes
g)\delta$.
\end{proposition}

\bp For $\lambda_1,\lambda_2,\lambda_3\in P$ we have to check that
$$
(\delta_{\lambda_1,\lambda_2}\otimes\iota)
\delta_{\lambda_1+\lambda_2,\lambda_3}=
(\iota\otimes\delta_{\lambda_2,\lambda_3})
\delta_{\lambda_1,\lambda_2+\lambda_3}.
$$
This reduces to showing that
$$
(m_{\mu_1,\mu_2,\lambda_1,\lambda_2}\otimes\iota
\otimes\iota)
m_{\mu_1+\mu_2,\mu_3,\lambda_1+\lambda_2,\lambda_3}
=(\iota\otimes\iota\otimes m_{\mu_2,\mu_3,\lambda_2,\lambda_3})
m_{\mu_1,\mu_2+\mu_3,\lambda_1,\lambda_2+\lambda_3},
$$
which follows immediately by definition.

\smallskip

Next we have to check that on $M_\lambda$ we have
$(\eps\otimes\iota)\delta_{0,\lambda}=\iota
=(\iota\otimes\eps)\delta_{\lambda,0}$. This is again straightforward.

\smallskip

Finally, to check commutativity of the diagram recall that the
isomorphism $$\Hom_{U_q\g}(M_{\lambda_1},V_\mu)\to V_\mu(\lambda_1)$$
comes from the homomorphisms $\Hom_{U_q\g}(\bar V_{\nu}\otimes
V_{\lambda_1+\nu},V_\mu)\to V_\mu(\lambda_1)$ given by $f\mapsto
f(\bar\xi_{\nu}\otimes\xi_{\lambda_1+\nu})$. It follows that it
suffices to check that
$$
m_{\nu,\omega,\lambda_1,\lambda_2}(\bar\xi_{\nu+\omega}\otimes
\xi_{\lambda_1+\lambda_2+\nu+\omega})=\bar\xi_{\nu}\otimes\xi_{\lambda_1
+\nu}\otimes\bar\xi_{\omega}\otimes\xi_{\lambda_2+\omega},
$$
but this is exactly the definition of $m$. \ep

The algebra $U_q\g$ acts by endomorphisms of the forgetful functor
$\CC(\g,\hbar)\to\Vect$. Our next goal is to show that the generators
of this action lift to endomorphisms of $M$.

Recall that in the previous section we defined morphisms
$$
\tau_{i;\eta,\mu}\colon V_{\eta+\mu-\alpha_i}\to V_\eta\otimes
V_\mu, \ \ \xi_{\eta+\mu-\alpha_i}\mapsto
[\eta(i)]_{q_i}\xi_\eta\otimes F_i\xi_\mu
-q_i^{\eta(i)}[\mu(i)]_{q_i}F_i\xi_\eta\otimes\xi_\mu.
$$
Similarly, define morphisms
$$
\bar\tau_{i;\mu,\eta}\colon\bar V_{\mu+\eta-\alpha_i}\to \bar V_\mu\otimes
\bar V_\eta, \ \ \bar\xi_{\mu+\eta-\alpha_i}\mapsto
[\eta(i)]_{q_i}E_i\bar\xi_\mu\otimes \bar\xi_\eta\
-q_i^{\eta(i)}[\mu(i)]_{q_i}\bar\xi_\mu\otimes E_i\bar\xi_\eta.
$$
Equivalently, $\bar\tau_{i;\mu,\eta}=\Sigma \tau_{i;\eta,\mu}$.

Consider the morphism
$$
\Psi^{\eta}_{i;\mu,\lambda+\alpha_i+\mu}
\colon \bar V_{\mu+\eta}\otimes V_{\lambda+\mu+\eta}
\to \bar V_{\mu}\otimes V_{\lambda+\alpha_i+\mu},\ \
\bar\xi_{\mu+\eta}\otimes\xi_{\lambda+\mu+\eta}
\mapsto \bar\xi_\mu\otimes F_i\xi_{\lambda+\alpha_i+\mu}.
$$
To see that it is well-defined, rewrite it as the composition
$$
\bar V_{\mu+\eta}\otimes V_{\lambda+\mu+\eta}
\xrightarrow{[\eta(i)]_{q_i}^{-1}\bar T_{\mu,\eta}\otimes
\tau_{i;\eta,\lambda+\alpha_i+\mu}} \bar V_{\mu}\otimes
\bar V_{\eta}\otimes V_\eta\otimes V_{\lambda+\alpha_i+\mu}
\xrightarrow{\iota\otimes S_\eta\otimes\iota}
\bar V_{\mu}\otimes V_{\lambda+\alpha_i+\mu}.
$$
Since $\bar\xi_\mu\otimes F_i\xi_{\lambda+\alpha_i+\mu}
=\Dhat_q(F_i)(\bar\xi_\mu\otimes\xi_{\lambda+\alpha_i+\mu})$, the
morphisms $\Psi$ are consistent with $\tr$ and hence define a morphism
$\tilde F_i\colon M_\lambda\to M_{\lambda+\alpha_i}$.

Similarly, consider the morphism
$$
\Phi^\eta_{i;\mu+\alpha_i,\lambda+\mu}\colon
\bar V_{\mu+\eta}\otimes V_{\lambda+\mu+\eta}
\to \bar V_{\mu+\alpha_i}\otimes
V_{\lambda+\mu},\ \
\bar\xi_{\mu+\eta}\otimes\xi_{\lambda+\mu+\eta}
\mapsto E_i\bar\xi_{\mu+\alpha_i}\otimes \xi_{\lambda+\mu},
$$
which can be equivalently written as the composition
$$
\bar V_{\mu+\eta}\otimes V_{\lambda+\mu+\eta}
\xrightarrow{[\eta(i)]_{q_i}^{-1}\bar\tau_{
i;\mu+\alpha_i,\eta}\otimes  T_{\eta,\lambda+\mu}}
\bar V_{\mu+\alpha_i}\otimes \bar V_{\eta}\otimes V_\eta\otimes
V_{\lambda+\mu} \xrightarrow{\iota\otimes
S_\eta\otimes\iota} \bar V_{\mu+\alpha_i}\otimes
V_{\lambda+\mu}.
$$
Again, using that $E_i\bar\xi_{\mu+\alpha_i}\otimes \xi_{\lambda+\mu}
=\Dhat_q(E_i)(\bar\xi_{\mu+\alpha_i}\otimes \xi_{\lambda+\mu})$, we see
that the morphisms $\Phi$ are consistent with $\tr$ and hence define a
morphism $\tilde E_i\colon M_\lambda\to M_{\lambda-\alpha_i}$.

Define also a morphism $\tilde K_i\colon M\to M$ by $\tilde
K_i|_{M_\lambda}=q_i^{\lambda(i)}$.

\begin{proposition} \label{pComonoid}
For all $1\le i\le r$ and $V\in\CC(\g,\hbar)$ the following diagrams
commute:
$$
\xymatrix{
\Hom_{U_q\g}(M,V) \ar[r]^{\mbox{\quad\quad\quad}\eta_V}\ar[d]_{\circ\tilde E_i}&
V \ar[d]^{E_i}\\
\Hom_{U_q\g}(M,V) \ar[r]_{\mbox{\quad\quad\quad}\eta_V} &
V}, \ \
\xymatrix{
\Hom_{U_q\g}(M,V) \ar[r]^{\mbox{\quad\quad\quad}\eta_V}\ar[d]_{\circ\tilde F_i}&
V \ar[d]^{F_i}\\
\Hom_{U_q\g}(M,V) \ar[r]_{\mbox{\quad\quad\quad}\eta_V} &
V},\ \
\xymatrix{
\Hom_{U_q\g}(M,V) \ar[r]^{\mbox{\quad\quad\quad}\eta_V}\ar[d]_{\circ\tilde K_i}&
V \ar[d]^{K_i}\\
\Hom_{U_q\g}(M,V) \ar[r]_{\mbox{\quad\quad\quad}\eta_V} &
V}.
$$
\end{proposition}

\bp To show commutativity of the first diagram it suffices to check
that if $$f\in\Hom_{U_q\g}(\bar V_{\mu+\alpha_i}\otimes
V_{\lambda+\mu},V)$$ then $ E_if(\bar \xi_{\mu+\alpha_i}\otimes
\xi_{\lambda+\mu}) =f\Phi^\eta_{i;\mu+\alpha_i,\lambda+\mu}(\bar
\xi_{\mu+\eta}\otimes \xi_{\lambda+\mu+\eta})$. Since
$$\Phi^\eta_{i;\mu+\alpha_i,\lambda+\mu}(\bar \xi_{\mu+\eta}\otimes
\xi_{\lambda+\mu+\eta})=\Dhat_q(E_i)(\bar \xi_{\mu+\alpha_i}\otimes
\xi_{\lambda+\mu}),$$ this is indeed true. The second diagram commutes
for similar reasons, while commutativity of the third diagram is
obvious. \ep

The next result will not be used later, but seems natural to complete
our discussion of the comonoid~$M$.

\begin{proposition}
There is a unital antihomomorphism $\pi\colon
U_q\g\mapsto\End_{U_q\g}(M)$ such that $\pi(X)=\tilde X$ for
$X\in\{E_i,F_i,K_i\}_{1\le i\le r}$. Furthermore, for any $\omega\in
U_q\g$ we have
$$
\delta\pi(\omega)=(\pi\otimes\pi)\Dhat_q(\omega)\delta\ \ \text{and}\ \
\eps\pi(\omega)=\hat\eps_q(\omega)\eps.
$$
\end{proposition}

\bp This is again a straightforward verification. Let us check for
example that
$$
\tilde F_j\tilde E_i-\tilde E_i\tilde F_j=\delta_{ij}\frac{\tilde K_i
-\tilde K_i^{-1}}{q_i-q_i^{-1}}.
$$
The vectors $\bar\xi_\mu\otimes\xi_{\lambda+\mu}$ define a
topologically cyclic vector $\Omega_\lambda\in M_\lambda$. The
morphisms $\tilde E_i$, $\tilde F_i$ and $\tilde K_i$ are defined by
$$
\tilde E_i\Omega_\lambda=E_i\Omega_{\lambda-\alpha_i},\ \
\tilde F_i\Omega_\lambda=F_i\Omega_{\lambda+\alpha_i}\ \ \text{and}\ \
\tilde K_i\Omega_\lambda=q_i^{\lambda(i)}\Omega_\lambda=K_i\Omega_\lambda.
$$
Therefore
$$
(\tilde F_j\tilde E_i-\tilde E_i\tilde F_j)\Omega_\lambda
=(E_iF_j-F_jE_i)\Omega_{\lambda-\alpha_i+\alpha_j}
=\delta_{ij}\frac{K_i -K_i^{-1}}{q_i-q_i^{-1}}\Omega_{\lambda-\alpha_i+\alpha_j}
=\delta_{ij}\frac{\tilde K_i -\tilde K_i^{-1}}{q_i-q_i^{-1}}
\Omega_{\lambda-\alpha_i+\alpha_j},
$$
which gives the result. \ep

\bigskip

\section{Proof of the main theorem}

We now return to the proof of Theorem~\ref{thmMain}. So let
$\E\in\U(G_q\times G_q)$ be a symmetric invariant $2$-cocycle
satisfying properties~\eqref{eHigh} and \eqref{eTauNorm}. Recall that
the latter properties mean that
$$
\E T_{\mu,\eta}=T_{\mu,\eta}\ \ \text{and}\ \
\E\tau_{i;\mu,\eta}=\tau_{i;\mu,\eta}.
$$
In the previous section we also introduced the maps $\bar T_{\mu,\eta}$
and $\bar\tau_{i;\mu,\eta}$. The first is an isomorphism of~$\bar
V_{\mu+\eta}$ onto the isotypic component of $\bar V_\mu\otimes \bar
V_\eta$ with lowest weight $-\mu-\eta$, that is, with highest
weight $\bar\mu+\bar\eta$. The second is an isomorphism of $\bar
V_{\mu+\eta-\alpha_i}$ onto the isotypic component with lowest
weight $-\mu-\eta+\alpha_i$, hence with highest weight
$\bar\mu+\bar\eta-\bar\alpha_i$. Therefore if we fix isomorphisms $\bar
V_\nu\cong V_{\bar \nu}$, then $\bar T_{\mu,\eta}$ and
$\bar\tau_{i;\mu,\eta}$ coincide with $T_{\bar\mu,\bar\eta}$ and
$\tau_{\bar i;\bar\mu,\bar\eta}$ up to scalar factors. Hence
properties~\eqref{eHigh} and \eqref{eTauNorm} also imply that
$$
\E \bar T_{\mu,\eta}=\bar T_{\mu,\eta}\ \ \text{and}\ \
\E\bar\tau_{i;\mu,\eta}=\bar\tau_{i;\mu,\eta}.
$$

Since $\E$ is invertible, the morphism $S_\mu\E\colon \bar
V_{\mu}\otimes V_\mu\to\C$ is nonzero, hence it is a nonzero multiple
of~$S_\mu$, so $S_\mu\E=\chi(\mu)S_\mu$ for some $\chi(\mu)\in\C^*$.
Explicitly, $\chi(\mu)=S_\mu\E(\bar\xi_{\mu}\otimes\xi_\mu)$.

\begin{lemma} \label{lCharact}
For all $\mu,\eta\in P_+$ and $\lambda\in P$ such that $\lambda+\mu\in
P_+$ we have
$\tr^\eta_{\mu,\lambda+\mu}\E=\chi(\eta)\E\tr^\eta_{\mu,\lambda+\mu}$.
\end{lemma}

\bp Applying $\iota\otimes\iota\otimes\Dhat_q$ to the cocycle identity
$$
(\E\otimes1)(\Dhat_q\otimes\iota)(\E)=(1\otimes\E)(\iota\otimes\Dhat_q)(\E),
$$
we get
$$
(\E\otimes1\otimes1)(\Dhat_q\otimes\Dhat_q)(\E)=(1\otimes(\iota\otimes\Dhat_q)(\E))
(\iota\otimes\Dhat_q^{(2)})(\E),
$$
where $\Dhat_q^{(2)}=(\iota\otimes\Dhat_q)\Dhat_q$. Replacing
$(\iota\otimes\Dhat_q)(\E)$ by
$(1\otimes\E^{-1})(\E\otimes1)(\Dhat_q\otimes\iota)(\E)$ on the right
hand side, we then get
\begin{equation*}
(\E\otimes\E)(\Dhat_q\otimes\Dhat_q)(\E)=
(1\otimes \E\otimes1)(1\otimes(\Dhat_q\otimes\iota)(\E))
(\iota\otimes\Dhat_q^{(2)})(\E),
\end{equation*}
which can also be written as
\begin{equation} \label{eCoDhat}
(\Dhat_q\otimes\Dhat_q)(\E)=
(1\otimes \E\otimes1)(1\otimes(\Dhat_q\otimes\iota)(\E))
(\iota\otimes\Dhat_q^{(2)})(\E)(\E^{-1}\otimes\E^{-1}),
\end{equation}
since $\E$ commutes with the image of $\Dhat_q$ by $G_q$-invariance.

We then compute
\begin{align*}
\tr^\eta_{\mu,\lambda+\mu}\E
&=(\iota\otimes S_\eta\otimes\iota)(\bar T_{\mu,\eta}\otimes
T_{\eta,\lambda+\mu})\E\\
&=(\iota\otimes S_\eta\otimes\iota)
(\Dhat_q\otimes\Dhat_q)(\E)
(\bar T_{\mu,\eta}\otimes T_{\eta,\lambda+\mu})\\
&=(\iota\otimes S_\eta\otimes\iota)
(1\otimes \E\otimes1)(1\otimes(\Dhat_q\otimes\iota)(\E))
(\iota\otimes\Dhat_q^{(2)})(\E)(\E^{-1}\otimes\E^{-1})
(\bar T_{\mu,\eta}\otimes T_{\eta,\lambda+\mu})\\
&=(\iota\otimes S_\eta\otimes\iota)
(1\otimes \E\otimes1)(1\otimes(\Dhat_q\otimes\iota)(\E))
(\iota\otimes\Dhat_q^{(2)})(\E)
(\bar T_{\mu,\eta}\otimes T_{\eta,\lambda+\mu})\\
&\quad\quad\quad (\hbox{by condition \eqref{eHigh}})\\
&=\chi(\eta)
(\iota\otimes S_\eta\otimes\iota)(\iota\otimes\Dhat_q^{(2)})(\E)
(\bar T_{\mu,\eta}\otimes T_{\eta,\lambda+\mu})\\
&\quad\quad\quad (\hbox{since}\ S_\eta\Dhat_q(\omega)=\hat\eps_q(\omega)S_\eta\
\hbox{and}\ (\hat\eps_q\otimes\iota)(\E)=1\ \ \hbox{by \eqref{eHigh}})\\
&=\chi(\eta)\E
(\iota\otimes S_\eta\otimes\iota)
(\bar T_{\mu,\eta}\otimes T_{\eta,\lambda+\mu})\\
&\quad\quad\quad (\hbox{since}\ S_\eta\Dhat_q(\omega)=\hat\eps_q(\omega)S_\eta\
\hbox{and}\ (\hat\eps_q\otimes\iota)\Dhat_q=\iota)\\
&=\chi(\eta)\E\tr^\eta_{\mu,\lambda+\mu}.
\end{align*} \ep

In particular, using that
$S_{\mu+\eta}=\tr^{\mu+\eta}_{0,0}=\tr^\mu_{0,0}\tr^\eta_{\mu,\mu}$ we
get
$$
\chi(\mu+\eta)S_{\mu+\eta}=S_{\mu+\eta}\E=S_\mu\tr^\eta_{\mu,\mu}\E
=\chi(\eta)S_\mu\E\tr^\eta_{\mu,\mu}
=\chi(\eta)\chi(\mu)S_\mu\tr^\eta_{\mu,\mu}
=\chi(\eta)\chi(\mu)S_{\mu+\eta}.
$$
Thus the map $\chi\colon P_+\to\C^*$ is a homomorphism, hence it
extends to a homomorphism $P\to\C^*$, which we continue to denote by
$\chi$. This together with the above lemma implies that the morphisms
$$\chi(\mu)^{-1}\E\colon \bar V_{\mu}\otimes V_{\lambda+\mu}\to
\bar V_{\mu}\otimes V_{\lambda+\mu}$$ are consistent with $\tr$, hence
define a morphism $\E_0\colon M_\lambda\to M_\lambda$. Note that $\E_0$
is invertible since $\E$ is.

\begin{lemma} \label{lCommut}
For all $1\le i\le r$ we have
$$
\tilde E_i\E_0=\chi(\alpha_i)\E_0\tilde E_i,\ \
\tilde F_i\E_0=\E_0\tilde F_i\ \ \text{and}\ \
\tilde K_i\E_0=\E_0\tilde K_i.
$$
\end{lemma}

\bp Recall that $\tilde E_i$ is defined using the morphisms
$\Phi^\eta_{i;\mu+\alpha_i,\lambda+\mu}$ given by the composition
$$
\bar V_{\mu+\eta}\otimes V_{\lambda+\mu+\eta}
\xrightarrow{[\eta(i)]_{q_i}^{-1}\bar\tau_{
i;\mu+\alpha_i,\eta}\otimes  T_{\eta,\lambda+\mu}}
\bar V_{\mu+\alpha_i}\otimes \bar V_{\eta}\otimes V_\eta\otimes
V_{\lambda+\mu} \xrightarrow{\iota\otimes
S_\eta\otimes\iota} \bar V_{\mu+\alpha_i}\otimes
V_{\lambda+\mu}.
$$
The same proof as that in Lemma~\ref{lCharact} shows that
$$
\Phi^\eta_{i;\mu+\alpha_i,\lambda+\mu}\E=\chi(\eta)
\E\Phi^\eta_{i;\mu+\alpha_i,\lambda+\mu}.
$$
The only difference is that $\bar T_{\mu,\eta}$ in that lemma gets
replaced by $\bar\tau_{ i;\mu+\alpha_i,\eta}$ and then instead of
condition~\eqref{eHigh} one uses condition~\eqref{eTauNorm}. Dividing
both sides of the above identity by $\chi(\mu+\eta)$, we get $\tilde
E_i\E_0=\chi(\alpha_i)\E_0\tilde E_i$.

Similarly, $\tilde F_i$ is defined using the morphisms
$\Psi^{\eta}_{i;\mu,\lambda+\alpha_i+\mu}$ given by the composition
$$
\bar V_{\mu+\eta}\otimes V_{\lambda+\mu+\eta}
\xrightarrow{[\eta(i)]_{q_i}^{-1}\bar T_{\mu,\eta}\otimes
\tau_{i;\eta,\lambda+\alpha_i+\mu}} \bar V_{\mu}\otimes
\bar V_{\eta}\otimes V_\eta\otimes V_{\lambda+\alpha_i+\mu}
\xrightarrow{\iota\otimes S_\eta\otimes\iota}
\bar V_{\mu}\otimes V_{\lambda+\alpha_i+\mu}.
$$
It follows that
$$
\Psi^{\eta}_{i;\mu,\lambda+\alpha_i+\mu}\E=\chi(\eta)
\E\Psi^{\eta}_{i;\mu,\lambda+\alpha_i+\mu},
$$
and dividing both sides by $\chi(\mu+\eta)$ we get $\tilde
F_i\E_0=\E_0\tilde F_i$.

The commutation with $\tilde K_i$ is obvious. \ep

The morphism $\E_0$ defines an endomorphism of the functor
$\Hom_{U_q\g}(M,\cdot)$. Since this functor is isomorphic to the
forgetful functor and the algebra of endomorphisms of the forgetful
functor is~$\U(G_q)$, the morphism $\E_0$ defines an invertible element
$c\in\U(G_q)$ such that for any $V\in\CC(\g,\hbar)$ the following
diagram commutes:
\begin{equation} \label{ec}
\xymatrix{
\Hom_{U_q\g}(M,V) \ar[r]^{\mbox{\quad\quad\quad}\eta_V}\ar[d]_{\circ\E_0}& V \ar[d]^c\\
\Hom_{U_q\g}(M,V) \ar[r]_{\mbox{\quad\quad\quad}\eta_V} &
V}.
\end{equation}

Proposition~\ref{pComonoid} and Lemma~\ref{lCommut} imply that
$$
c E_i=\chi(\alpha_i)E_i c,\ \ c F_i=F_ic\ \ \text{and}\ \ cK_i=K_ic.
$$
Since $E_iF_i-F_iE_i$ coincides with $K_i-K_i^{-1}$ up to a scalar
factor, this is possible only if $\chi(\alpha_i)=1$. Therefore $c$
belongs to the center of $\U(G_q)$.

\begin{lemma} \label{lcocycle}
We have $\delta\E_0=\E(\E_0\otimes\E_0)\delta$.
\end{lemma}

\bp Recall that $\delta$ is defined using the morphisms
$$
m_{\mu,\eta,\lambda_1,\lambda_2}\colon \bar V_{\mu+\eta}\otimes
V_{\lambda_1+\lambda_2+\mu+\eta}\to \bar V_{\mu}\otimes
V_{\lambda_1+\mu}\otimes \bar V_{\eta}\otimes V_{\lambda_2+\eta}
$$
given by
$m_{\mu,\eta,\lambda_1,\lambda_2}=q^{(\lambda_1+\mu,\eta)}(\iota\otimes
\sigma \otimes\iota)(\bar T_{\mu,\eta}\otimes
T_{\lambda_1+\mu,\lambda_2+\eta})$. The same computation as that in
Lemma~\ref{lCharact} shows that
\begin{align*}
m_{\mu,\eta,\lambda_1,\lambda_2}\E
&=q^{(\lambda_1+\mu,\eta)}(\iota\otimes
\sigma \otimes\iota)(1\otimes \E\otimes1)(1\otimes(\Dhat_q\otimes\iota)(\E))
(\iota\otimes\Dhat_q^{(2)})(\E)(\bar T_{\mu,\eta}\otimes
T_{\lambda_1+\mu,\lambda_2+\eta})\\
&=q^{(\lambda_1+\mu,\eta)}(1\otimes \E\otimes1)(1\otimes(\Dhat_q\otimes\iota)(\E))
(\iota\otimes\Dhat_q^{(2)})(\E)(\iota\otimes
\sigma \otimes\iota)(\bar T_{\mu,\eta}\otimes
T_{\lambda_1+\mu,\lambda_2+\eta})\\
&\quad\quad\quad(\hbox{since}\ \sigma\E=\E\sigma\ \hbox{by the assumption that}
\ \E\ \hbox{is symmetric})\\
&=q^{(\lambda_1+\mu,\eta)}(\Dhat_q\otimes\Dhat_q)(\E)(\E\otimes\E)(\iota\otimes
\sigma \otimes\iota)(\bar T_{\mu,\eta}\otimes
T_{\lambda_1+\mu,\lambda_2+\eta})\\
&\quad\quad\quad(\hbox{by}\ \eqref{eCoDhat})\\
&=(\Dhat_q\otimes\Dhat_q)(\E)(\E\otimes\E)m_{\mu,\eta,\lambda_1,\lambda_2},
\end{align*}
which proves the lemma. \ep

\bp[Proof of Theorem~\ref{thmMain}] For $U,V\in\CC(\g,\hbar)$ and
$f\in\Hom_{U_q\g}(M,U)$, $g\in\Hom_{U_q\g}(M,V)$ we have
\begin{align*}
\Dhat_q(c)(\eta_U(f)\otimes \eta_V(g))&=\Dhat_q(c)\eta_{U\otimes V}((f\otimes g)\delta)
&&(\hbox{by Proposition~\ref{pComonoid0}})\\
&=\eta_{U\otimes V}((f\otimes g)\delta\E_0) &&(\hbox{by}\ \eqref{ec})\\
&=\eta_{U\otimes V}((f\otimes g)\E(\E_0\otimes\E_0)\delta)
&&(\hbox{by Lemma~\ref{lcocycle}})\\
&=\eta_{U\otimes V}(\E(f\otimes g)(\E_0\otimes\E_0)\delta)\\
&=\E\eta_{U\otimes V}((f\E_0\otimes g\E_0)\delta)
&&(\hbox{by naturality of}\ \eta\ \hbox{and}\ G_q\hbox{-invariance of}\ \E)\\
&=\E(\eta_U(f\E_0)\otimes \eta_V(g\E_0))
&&(\hbox{by Proposition~\ref{pComonoid0}})\\
&=\E(c\otimes c)(\eta_U(f)\otimes \eta_V(g))
&&(\hbox{by}\ \eqref{ec}).
\end{align*}
It follows that $\Dhat_q(c)=\E(c\otimes c)$. \ep

With a bit more work one can show that in fact $\E=1$.

\begin{corollary} \label{cE1}
If $\E\in\U(G_q\times G_q)$ is a symmetric invariant $2$-cocycle
satisfying properties \eqref{eHigh} and \eqref{eTauNorm}, then $\E=1$.
\end{corollary}

\bp By Theorem~\ref{thmMain}, $\E$ is the coboundary of a central
element, that is, $\E=(c\otimes c)\Dhat_q(c)^{-1}$. The element $c$
acts on $V_\mu$ by a scalar $\chi(\mu)$. Condition~\eqref{eHigh} means
then that $\chi\colon P_+\to\C^*$ is a homomorphism, so $\chi$ extends
to a homomorphism $P\to\C^*$. Condition~\eqref{eTauNorm} implies that
$\chi(\alpha_i)=1$ for all $i$, so $\chi$ is trivial on the root
lattice $Q$. In other words, $\chi$ is a character of $P/Q$. But then
$c$ is group-like, that is, $\Dhat_q(c)=c\otimes c$. This is well-known
for $q=1$, since the characters of $P/Q$ are in a one-to-one
correspondence with the elements of the center of $G$, see
e.g.~\cite[Theorem~26.3]{Bu}, and so $c$ belongs to $G\subset\U(G)$ and
is therefore group-like. For general $q$, the canonical identification
of the centers of~$\U(G_q)$ and~$\U(G)$ extends to an isomorphism of
algebras, since the dimensions of the irreducible representations with
a given highest weight do not depend on $q$. Since the fusion rules do
not depend on $q$ either, there exists $\F\in\U(G\times G)$ such that $
\Dhat_q=\F\Dhat(\cdot)\F^{-1}$. Then as $c$ is group-like in
$(\U(G),\Dhat)$, we have
$$
\Dhat_q(c)=\F\Dhat(c)\F^{-1}=\F(c\otimes c)\F^{-1}=c\otimes c,
$$
and so $c$ is group-like in $(\U(G_q),\Dhat_q)$ as well. Hence
$\E=(c\otimes c)\Dhat_q(c)^{-1}=1$. \ep

In the above proof we remarked that a central element in $\U(G_q)$ is
group-like if it is defined by a character of $P/Q$. The converse is
also true.

\begin{proposition}
A central element of $\U(G_q)$ is group-like if and only if it is
defined by a character of $P/Q$.
\end{proposition}

\bp We only have to show that if $c$ is central and group-like then it
is defined by a character of~$P/Q$. Similarly to the proof of
Corollary~\ref{cE1}, we first conclude that $c$ is group-like in
$\U(G)$ as well. Then~$c$ is a central element of the complexification
$G_\C\subset\U(G)$ of $G$ by~Theorem~\ref{thmGroupLike}. Hence it
belongs to~$G$ and is defined by a character of~$P/Q$, see
again~\cite[Theorem~26.3]{Bu}. \ep

\begin{corollary}
\label{coruniqunec}
If $\E\in\U(G_q\times G_q)$ is a symmetric invariant $2$-cocycle then a
central element $c$ such that $\E=(c\otimes c)\Dhat_q(c)^{-1}$ is
defined uniquely up to a character of $P/Q$.
\end{corollary}

\bigskip

\section{Uniqueness of the Drinfeld twist}
\label{SecTwist}

We shall assume throughout this section that $q>0$. Let $\hbar\in i\R$
be such that $q=e^{\pi i\hbar}$. Denote by $t\in\g\otimes\g$ the
$\g$-invariant symmetric element defined by the Killing form normalized
so that the induced form on $\h^*$ satisfies $(\alpha, \alpha)=2$ for
short roots. Let $\Phi_{KZ}=\Phi(\hbar t_{12},\hbar
t_{23})\in\U(G\times G\times G)$ be the Drinfeld associator defined via
monodromy of the KZ-equations, see e.g.~\cite{NT3} for details.

Recall that from the work of Kazhdan and Lusztig~\cite{KL1} one can
derive~\cite{NT3} the following analytic version of a famous result of
Drinfeld~\cite{Dr1,Dr2}.

\begin{theorem} \label{Dr}
For any isomorphism $\varphi\colon\U(G_q)\to\U(G)$ extending the
canonical identification of the centers there exists an invertible
element $\F\in\U(G\times G)$ such that \enu{i}
$(\varphi\otimes\varphi)\Dhat_q=\F\Dhat\varphi(\cdot)\F^{-1}$; \enu{ii}
$(\hat\eps\otimes\iota)(\F)=(\iota\otimes\hat\eps)(\F)=1$; \enu{iii}
$(\varphi\otimes\varphi)(\RR_\hbar)=\F_{21}q^t\F^{-1}$; \enu{iv}
$\Phi_{KZ}=(\iota\otimes\Dhat)(\F^{-1})
(1\otimes\F^{-1})(\F\otimes1)(\Dhat\otimes\iota)(\F)$.

\medskip

\noindent In addition, if $\varphi$ is a $*$-isomorphism then $\F$ can
be taken to be unitary.
\end{theorem}

Any such element $\F$ is called a Drinfeld twist. Our next result
asserts that for $\varphi$ fixed, the Drinfeld twist is unique up to
coboundary of a central element. This is an equivalent form of Theorem
\ref{thmMain}.

\begin{theorem}
\label{corunique} Suppose $\F$ and $\F'$ are two Drinfeld twists for
the same isomorphism $\varphi$. Then there exists a central element $c$
of $\U(G)$ such that $\F' =(c\otimes c)\F\Dhat(c)^{-1}$. When $\varphi$
is a $*$-isomorphism and both Drinfeld twists are unitary, then $c$ can
also be chosen to be unitary.
\end{theorem}

\bp To simplify the notation we shall omit $\varphi$ in the
computations, so we identify $\U(G_q)$ and~$\U(G)$ as algebras. Set
$\E=\F'\F^{-1}$. Then
$$
(\iota\otimes\Dhat)(\F^{-1})
(1\otimes\F^{-1})(\F\otimes1)(\Dhat\otimes\iota)(\F)
=(\iota\otimes\Dhat)(\F^{-1}\E^{-1})
(1\otimes\F^{-1}\E^{-1})(\E\F\otimes1)(\Dhat\otimes\iota)(\E\F).
$$
Multiplying by $(1\otimes\F)(\iota\otimes\Dhat)(\F)$ on the left and by
$(\Dhat\otimes\iota)(\F^{-1})(\F^{-1}\otimes1)$ on the right, and using
that $\F\Dhat(\cdot)\F^{-1}=\Dhat_q$, we get
$$
1=(\iota\otimes\Dhat_q)(\E^{-1})
(1\otimes\E^{-1})(\E\otimes1)(\Dhat_q\otimes\iota)(\E).
$$
Therefore $\E$ is a $2$-cocycle for $(\U(G_q),\Dhat_q)$. Since
$$
\E\Dhat_q(\cdot)\E^{-1}=\E\F\Dhat(\cdot)\F^{-1}\E^{-1}
=\F'\Dhat(\cdot)\F'^{-1}=\Dhat_q,
$$
the cocycle $\E\in\U(G_q\times G_q)$ is invariant, and since
$$
\E_{21}\RR_\hbar\E^{-1}=\E_{21}\F_{21}q^t\F^{-1}\E^{-1}
=\F'_{21}q^t\F'^{-1}=\RR_\hbar,
$$
it is symmetric. By Theorem~\ref{thmMain} there exists a central
element $c$ of $\U(G_q)=\U(G)$ such that $$\E=(c\otimes
c)\Dhat_q(c)^{-1},$$ so that $ \F'=(c\otimes
c)\Dhat_q(c^{-1})\F=(c\otimes c)\F\Dhat(c^{-1}) $, and the first claim
is proved. The second claim is immediate from Lemma \ref{lempolar}. \ep

As for dependence on $\varphi$, if $\varphi'\colon\U(G_q)\to\U(G)$ is
another isomorphism extending the canonical identification of the
centers, there exists an invertible element $u$ of $\U(G)$ such that
$\varphi' =u\varphi(\cdot)u^{-1}$, and then $\F_u =(u\otimes u)\F\Dhat
(u)^{-1}$ is a Drinfeld twist for $\varphi'$. By
Theorem~\ref{corunique} all Drinfeld twists for~$\varphi'$ will
therefore be cohomologous to $\F_u$. So up to coboundary, there is only
one Drinfeld twist irrespectively of the choice of the isomorphism
$\varphi$. When one considers a $*$-isomorphism $\varphi$ together with
a unitary Drinfeld twist $\F$, then $u$ can be chosen to be unitary,
and consequently, irrespectively of $\varphi$, there is only one
unitary Drinfeld twist up to coboundary of a unitary element.



\medskip

In the language of cohomology from Section \ref{seccohom}, the Drinfeld
associator $\Phi=\Phi_{KZ}$ is a unitary counital invariant
$3$-cocycle for $(\U(G),\Dhat)$ satisfying the 
equation
$$\RR_{12}\Phi_{312}\RR_{13}\Phi^{-1}_{132}\RR_{23}\Phi_{123}
=\Phi_{321}\RR_{23}\Phi^{-1}_{231}\RR_{13}\Phi_{213}\RR_{12},$$ which
is some sort of symmetry condition. Theorem \ref{Dr} tells us then that
$\Phi =\partial (\F^{-1})$, where the coboundary operator $\partial$
refers to $(\U(G),\F\Dhat (\cdot)\F^{-1})$, which is isomorphic to
$(\U(G_q), \Dhat_q)$. This should be compared with Theorem
\ref{thmMain} stating that any symmetric invariant $2$-cocycle for
$(\U(G_q),\Dhat_q)$ is the coboundary of a central element.

\bigskip

\section{Uniqueness of the Dirac operator}\label{sdirac}

As in the previous section, we assume that $q>0$ and $\hbar\in i\R$ is
such that $q=e^{\pi i\hbar}$. In \cite{NT2} we constructed a quantum
Dirac operator $D_q$ on $G_q$ that defines a biequivariant spectral
triple which is an isospectral deformation of that defined by the Dirac
operator $D$ on $G$. We briefly recall this construction.

The Riemannian metric on $G$ is defined using the invariant form
$-(\cdot,\cdot)$ on $\g$. Consider a basis~$\{x_i\}_i$ of $\g$ such
that $(x_i ,x_j)=-\delta_{ij}$, and let $\gamma\colon\g\to\Clg$ denote
the inclusion of $\g$ into the complex Clifford algebra with the
convention that $\gamma (x_i)^2=-1$. Identifying ${\mathfrak{so} }(\g)$
with ${\mathfrak{spin}}(\g)$, the adjoint action is defined by the
representation $\add\colon\g\to{\mathfrak{spin}} (\g)\subset\Clg$ given
by
$$\add(x)=\frac{1}{4}\sum_i \gamma (x_i )\gamma ([x,x_i ]).$$
We denote by the same symbol $\add$ the corresponding homomorphism
$\U(G)\to\Clg$.

Let $s\colon\Clg\to \End(\Sp)$ be an irreducible representation. Denote
by $\partial$ the representation of $U\g$ by left-invariant
differential operators. Identifying the sections $\Gamma (S)$ of the
spin bundle $S$ over $G$ with $C^\infty (G)\otimes\Sp$, the Dirac
operator $D\colon C^\infty (G)\otimes\Sp\to C^\infty (G)\otimes\Sp$
defined using the Levi-Civita connection, can be written as
$D=(\partial\otimes s)(\D)$, where $\D\in U\g\otimes\Clg$ is given by
the formula
$$\D=\sum_i (x_i\otimes\gamma (x_i )+\frac{1}{2}\otimes\gamma(x_i)\add(x_i)).$$


\smallskip

Denote by $\C[G_q]$ the linear span of matrix coefficients of finite
dimensional admissible representations of $U_q\g$. It is a Hopf
$*$-algebra with comultiplication $\Delta_q$, and $\U(G_q)$ is its dual
space. Let $(L^2 (G_q),\pi_{r,q},\xi_q )$ be the GNS-triple defined by
the Haar state on $\C [G_q ]$. The left and right regular
representations of~$W^*(G_q)$ on $L^2(G_q)$, denoted by $\hat\pi_{r,q}$
and $\partial_q$ correspondingly, are defined by
$$
\hat\pi_{r,q}(\omega)\pi_{r,q}(a)\xi_q=(\omega
S^{-1}_q\otimes\pi_{r,q})\Delta_q(a)\xi_q,
$$
where $S_q$ is the antipode on $\C[G_q]$, and
\begin{equation*} \label{eRegL}
\partial_q(\omega)\pi_{r,q}(a)\xi_q
=(\pi_{r,q}\otimes\omega)\Delta_q(a)\xi_q
=a_{(1)}(\omega)\pi_{r,q}(a_{(0)})\xi_q.
\end{equation*}

Pick a $*$-isomorphism $\varphi\colon\U(G_q)\to\U(G)$ and a unitary
Drinfeld twist $\F\in\U(G\times G)$.
The quantum Dirac operator $D_q$ is the unbounded operator on
$L^2(G_q)\otimes\Sp$ defined by
$$
D_q=(\partial_q\otimes s)(\D_q),
$$
where
$\D_q\in\U(G_q)\otimes\Clg$ is given by
$$
\D_q=(\varphi^{-1}\otimes\iota)((\iota\otimes\add)(\F)
\D(\iota\otimes\add)(\F^*)).
$$
The operator $D_q$ is $G_q$-biequivariant in the sense that it commutes
with all operators of the form $\hat\pi_{r,q}(x)\otimes1$ and
$(\partial_q\times s\,\add_q)(x)$, $x\in W^*(G_q)$, where
$\add_q=\add\varphi$.


\begin{theorem}\label{coruniqueD}\mbox{\ }
\enu{i} For fixed $\varphi$ the Dirac operator $D_q$ does not depend on
the chosen Drinfeld twist $\F$. \enu{ii} The biequivariant spectral
triple $(\C[G_q],L^2(G_q)\otimes\Sp,D_q)$ does not depend on the choice
of $\varphi$ and~$\F$ up to unitary equivalence.
\end{theorem}

\bp By Theorem \ref{corunique} any other unitary Drinfeld twist $\tilde
\F$ for the same $\varphi$ has the form $$\tilde\F=(c\otimes c)\F\Dhat
(c)^*$$ for a central unitary element $c$ of $\U(G)$. Denoting the
element $\D_q$ defined by $\tilde\F$ by $\tilde\D_q$, we get
\begin{align*}
(\varphi\otimes\iota)(\tilde\D_q)
&=(\iota\otimes\add)((c\otimes c)\F\Dhat(c)^*)\D(\iota\otimes\add)(\Dhat(c)
\F^*(c^*\otimes c^*))& & \\
&=(\iota\otimes\add)((c\otimes c)\F)\D(\iota\otimes\add)(\F^*(c^*\otimes c^*))&
&(\text{since}\ \D\ \text{is}\ \g\text{-invariant})\\
&=(1\otimes\add(c))(\iota\otimes\add)(\F)\D(\iota\otimes\add)(\F^*)
(1\otimes\add(c^*))& &(\text{since}\ c\ \text{is central})\\
&=(\iota\otimes\add)(\F)\D(\iota\otimes\add)(\F^*)
& &(\text{since}\ \add(c)\ \text{is a scalar}),
\end{align*}
where in the last step we used the known fact that $\add$ is a multiple
of an irreducible representation, namely, of the representation with
highest weight $\rho$, half the sum of positive roots. This proves (i).

\smallskip

If we choose another $*$-isomorphism $\varphi'\colon\U(G_q)\to\U(G)$,
there exists a unitary $u$ such that $\varphi'=u\varphi(\cdot)u^*$. We
can take the element $\F'=(u\otimes u)\F\Dhat(u^*)$ as a unitary
Drinfeld twist for $\varphi'$. Then for the element $\D_q'$ defined by
$\varphi'$ and $\F'$ we get, using that $\D$ commutes with
$(\iota\otimes\add)\Dhat(u)$, that
$$
\D_q'=(\varphi'^{-1}\otimes\iota)\left(
(\iota\otimes\add)((u\otimes u)\F)\D(\iota\otimes\add)(\F^*(u^*\otimes u^*))\right)
=(1\otimes\add(u))\D_q(1\otimes\add(u^*)).
$$
We also have $\add'_q:=\add\varphi'=\add(u)\add_q(\cdot)\add(u^*)$.
Therefore the operator $1\otimes s\, \add(u)$ provides a unitary
equivalence between the biequivariant spectral triples
$(\C[G_q],L^2(G_q)\otimes\Sp,D_q)$ and
$(\C[G_q],L^2(G_q)\otimes\Sp,D_q')$. \ep

\bigskip

\appendix

\section{}

Let $G$ be a compact Lie group and $G_\C$ its analytic
complexification. By definition, any continuous finite dimensional
representation $G\to\GL(V)$ extends uniquely to a holomorphic
representation $G_\C\to\GL(V)$. Hence every element $g\in G_\C$ can be
considered as an element of~$\U(G)$. Furthermore, $\Dhat(g)=g\otimes g$
by analyticity, since this is true for all $g\in G$ (to be more
precise, in order to not worry about topology, we should first apply a
finite dimensional representation $\pi_1\otimes\pi_2$).
Therefore~$G_\C$ consists of group-like elements. We will show
that these are all; for $G=\SU(n)$ this is~\cite[Theorem~2]{Bi}.

\begin{theorem} \label{thmGroupLike}
For any compact Lie group $G$ the set of group-like elements in $\U(G)$
coincides with~$G_\C$.
\end{theorem}

\bp Let $a\in\U(G)$ be a group-like element, so $a$ is invertible and
$\Dhat(a)=a\otimes a$. Assume first that~$a$ is bounded, so it belongs
to the von Neumann algebra $W^*(G)\subset\U(G)$ of $G$. Then $a\in
G\subset G_\C$. This is a well-known result going back to
Tatsuuma~\cite{Tat} and valid for any locally compact group. Here is a
short proof.

Consider $W^*(G)$ as the von Neumann algebra generated by the operators $\lambda_g$ of the left regular representation of $G$. Therefore we want to prove that $a=\lambda_g$ for some $g\in G$. Let $U$ be an open
neighbourhood of the unit element $e\in G$. Consider the set $K_U$
consisting of all elements $g\in G$ for which there exists a function
$f\in L^2(G)$ with essential support in $U$ such that
$(\operatorname{ess.supp} af)\cap gU\ne\emptyset$. As the whole space
$L^2(G)$ is spanned by right translations of functions with essential
support in $U$, there exists $f$ with $\operatorname{ess.supp} f\subset
U$ such that $af\ne0$. It follows that $K_U$ is non-empty. We claim
that if $g_0\in K_U$ then
\begin{itemize}
\item[(i)] the element $a$ lies in the strong operator closure of
    the span of $\lambda_g$ with $g\in g_0UU^{-1}$;
\item[(ii)] $K_U\subset\overline{g_0UU^{-1}UU^{-1}}$.
\end{itemize}
Indeed, consider functions $f$ and $h$ such that
$\operatorname{ess.supp} f\subset U$, $\operatorname{ess.supp} h\subset
g_0U$ and $(af,h)\ne0$. Denote by $\omega$ the normal linear functional
$(\cdot\,f,h)$ on $B(L^2(G))$. Then
$$
(\iota\otimes\omega)\Dhat(a)=(\iota\otimes\omega)(a\otimes a)=\omega(a)a,
$$
and on the other hand,
$$
(\iota\otimes\omega)\Dhat(\lambda_g)=\omega(\lambda_g)\lambda_g=0\ \
\text{for}\ \ g\notin g_0UU^{-1}.
$$
Since $a$ can be approximated by linear combinations of the operators
$\lambda_g$, applying the normal operator $(\iota\otimes\omega)\Dhat$
to these approximations we get (i). Now if $f\in L^2(G)$ is arbitrary
with $\operatorname{ess.supp} f\subset\nolinebreak U$, by (i) we have
$\operatorname{ess.supp} af\subset \overline{g_0UU^{-1}U}$, whence
$K_U\subset \overline{g_0UU^{-1}UU^{-1}}$.

If $V\subset U$ are two neighbourhoods of $e\in G$ then clearly
$K_V\subset K_U$. Property (ii) implies that the intersection of the
sets $\overline K_U$ consists of exactly one point, which we denote by
$g_0$. Property (i) implies that $a$ belongs to the strong operator
closure of the span of the operators $\lambda_g$ with $g$ lying in an
arbitrarily small neighbourhood of $g_0$. We want to prove that this
forces $a=\lambda_{g_0}$. Replacing $a$ by~$\lambda^{-1}_{g_0}a$ we may
assume that $g_0=e$. Then for any $f\in L^2(G)$ we get
$$
\operatorname{ess.supp} af\subset \operatorname{ess.supp} f.
$$
If we consider the action of $L^\infty(G)$ on $L^2(G)$ by
multiplication, by regularity of the Haar measure this implies that $a$
commutes with the characteristic function of any measurable set. It
follows that $a\in L^\infty(G)$. Since $a$ commutes with the operators
of the right regular representation, this implies that $a$ is a scalar,
and since it is group-like, we get $a=1$.

\smallskip

Consider now an arbitrary group-like element $a\in\U(G)$. Then $a^*a$
is group-like as well, hence $|a|$ is also group-like. It follows that
if $a=u|a|$ is the polar decomposition then $u$ is group-like. By the
first part of the proof we know that $u\in G$. So we just have to show
that $|a|\in G_\C$. In other words, we may assume that $a$ is positive.

For every $z\in\C$ we have
$$
\Dhat(a^z)=\Dhat(a)^z=(a\otimes a)^z=a^z\otimes a^z.
$$
In particular, the bounded elements~$a^{it}$, $t\in\R$, are group-like,
hence they lie in $G\subset\U(G)$. It follows that there exists
$X\in\g$ such that $a^{it}=\exp{tX}$ for~$t\in\R$, whence
$a^z=\exp({-izX})\in G_\C$ for all $z\in\C$, since both $a^z$ and
$\exp({-izX})$ are analytic functions in $z$ which coincide for
$z\in i\R$. In particular, $a=\exp({-iX})\in G_\C$. \ep

\bigskip

\section{}

The proof of the main theorem can also be applied in the formal
deformation setting. However, in this case it is easier to follow
Drinfeld's cohomological arguments for $3$-cocycles~\cite{Dr2}, see
also the proof of~\cite[Theorem~XVIII.8.1]{Ka}. Although a translation
of those arguments into our setting of $2$-cocycles is completely
straightforward, we include it in this appendix for the reader's
convenience.

\begin{theorem}
Let $\g$ be a finite dimensional semisimple Lie algebra, and let
$(U\g[[h]],\Dhat_h,\RR_h)$ be a quasitriangular deformation of
$(U\g,\Dhat)$. Assume $\E\in(U\g\otimes U\g)[[h]]$ is a symmetric
invariant $2$-cocycle such that $\E=1\mod h$, so \enu{i}
$[\E,\Dhat_h(a)]=0$ for all $a\in U\g[[h]]$; \enu{ii}
$\RR_h\E=\E_{21}\RR_h$; \enu{iii}
$(\E\otimes1)(\Dhat_h\otimes\iota)(\E)=(1\otimes\E)(\iota\otimes\Dhat_h)(\E)$.

\smallskip\noindent
Then there exists a central element $c\in U\g[[h]]$ such that $c=1 \mod
h$ and $\E=(c\otimes c)\Dhat_h(c)^{-1}$.
\end{theorem}

\bp We will construct by induction central elements $c_n\in U\g[[h]]$,
$n\ge0$, such that $c_0=1$ and
$$
\E=(c_n\otimes c_n)\Dhat_h(c_n)^{-1}\mod h^{n+1}\ \ \text{and}\ \
c_n=c_{n-1}\mod h^n\ \ \text{for}\ \ n\ge1.
$$
Then the sequence $\{c_n\}_n$ converges to the required element $c$.

Assume $c_0,\dots,c_{n-1}$ are constructed. Let $\varphi\in U\g\otimes
U\g$ be such that
$$
\E=(c_{n-1}\otimes c_{n-1})\Dhat_h(c_{n-1})^{-1}+h^n\varphi\mod h^{n+1}.
$$
Reducing conditions (i)-(iii) modulo $h^{n+1}$ and using that
$\Dhat_h=\Dhat$, $\RR_h=1$ and $c_{n-1}=1$ modulo $h$, we get
$$
[\varphi,\Dhat(a)]=0 \ \text{for}\ a\in U\g,\ \ \varphi=\varphi_{21}\ \
\text{and}\ \
\varphi\otimes1+(\Dhat\otimes\iota)(\varphi)=1\otimes \varphi+
(\iota\otimes\Dhat)(\varphi).
$$
In the notation of \cite[Ch.~XVIII.5]{Ka} the last two identities mean
that $\varphi$ is a $2$-cocycle in the complex $(T_-(U\g),\delta)$.
Since the symmetrization map $\eta\colon S\g\to U\g$ is an isomorphism
of coalgebras, we have $H^{2k}(T_-(U\g),\delta)=0$ for all $k\ge0$ by
\cite[Theorem~XVIII.7.1]{Ka}. Therefore $\varphi$ is the coboundary of
an element $f\in U\g$, so that
$$
\varphi=f\otimes1+1\otimes f-\Dhat(f).
$$
Furthermore, since $\varphi$ is $\g$-invariant, by
\cite[Proposition~XVIII.6.2]{Ka} we can choose $f$ to be $\g$-invariant
as well. Then we put $c_n=(1+h^nf)c_{n-1}$. \ep

\end{document}